\documentclass[12pt,a4paper]{amsart}
\usepackage{pstricks,pstcol,pst-plot,pst-node,pst-tree,amssymb}
\usepackage{pstricks-add}

\newtheorem{theorem}{Theorem}

\DeclareMathOperator{\Imaginary}{Im}
\DeclareMathOperator{\Real}{Re}

\begin{document}
\title[Q-prime curvature]{An integral formula for the $Q$-prime curvature in
3-dimensional CR geometry}
\author{Jeffrey S. Case}
\address{Department of Mathematics, \\
Pennsylvania State University \\
University Park, PA 16802, U.S.A.}
\email{jqc5026@psu.edu}
\thanks{}
\author{Jih-Hsin Cheng}
\address{Institute of Mathematics, Academia Sinica, Taipei and National
Center for Theoretical Sciences, Taipei Office, Taiwan, R.O.C.}
\email{cheng@math.sinica.edu.tw}
\thanks{}
\author{Paul Yang}
\address{Department of Mathematics, Princeton University, Princeton, NJ
08544, U.S.A.}
\email{yang@Math.Princeton.EDU}
\urladdr{}
\subjclass[2010]{Primary 32V05; Secondary 32V20}
\keywords{$Q$-prime curvature; pseudo-Einstein manifold; CR manifold}
\thanks{}

\begin{abstract}
 We give an integral formula for the total $Q^\prime$-curvature of a three-dimensional CR manifold with positive CR Yamabe constant and nonnegative Paneitz operator.  Our derivation includes a relationship between the Green's functions of the CR Laplacian and the $P^\prime$-operator.
\end{abstract}

\maketitle
\date{}


\section{Introduction}

The $Q^\prime$-curvature, introduced to three-dimensional CR manifolds by the first- and third-named authors~\cite{CY} and to higher-dimensional CR manifolds by Hirachi~\cite{H2}, has recently emerged as the natural CR counterpart to Branson's $Q$-curvature in conformal geometry.  The analogies are especially strong in dimension three, where it is known that the total $Q^\prime$-curvature is a biholomorphic invariant --- indeed, it is a multiple of the Burns--Epstein invariant~\cite{CY,CL} --- and gives rise to a CR invariant characterization of the standard CR three-sphere.

The above discussion is complicated by the fact that the $Q^\prime$-curvature is most naturally defined only for pseudo-Einstein contact forms.  A pseudohermitian manifold $(M^3,J,\theta)$ is \emph{pseudo-Einstein} if the curvature $R$ and torsion $A_{11}$ of the Tanaka--Webster connection satisfy $R_1=iA_{11,\bar 1}$.  It is known~\cite{H1} that if $\theta$ is a pseudo-Einstein contact form, then $\hat\theta:=e^\Upsilon\theta$ is pseudo-Einstein if and only if $\Upsilon$ is a CR pluriharmonic function.  Moreover, if $M^3$ is embedded in $\mathbb{C}^2$, then pseudo-Einstein contact forms arise from solutions of Fefferman's Monge--Amp\`ere equation~\cite{FH}.  For a pseudo-Einstein manifold $(M^3,J,\theta)$, the $Q^\prime$-curvature is defined by
\[ Q^\prime := -2\Delta_b R + R^2 - 4\lvert A_{11}\rvert^2 . \]
The behavior of $Q^\prime$ under the conformal transformation of $\theta$ to $\hat\theta$ is controlled by the $P$-prime operator $P^\prime$ and the Paneitz operator $P$, which have the local expressions
\begin{align*}
 P^\prime(u) & := 4\Delta_b^2 - 8\Imaginary (A_{11}u_{\bar 1})_{;\bar1} - 4\Real (Ru_1)_{;\bar 1}, \\
 P(u) & := \Delta_b^2u + T^2u - 4\Imaginary (A_{11}u_{\bar 1})_{;\bar1} .
\end{align*}
More precisely, if $\hat\theta=e^\Upsilon\theta$ and $\theta$ are both pseudo-Einstein, then
\begin{equation}
 \label{eqn:qprime_transformation}
 e^{2\Upsilon}\hat Q^\prime = Q^\prime + P^\prime(\Upsilon) + \frac{1}{2}P(\Upsilon^2) .
\end{equation}
From this formula, it is clear that the total $Q^\prime$-curvature is independent of the choice of pseudo-Einstein contact form.  A direct computation also shows that if the holomorphic tangent bundle of $M$ is trivial, then the total $Q^\prime$-curvature is a multiple of the Burns--Epstein invariant~\cite{CL}.

The CR Yamabe constant of a CR manifold $(M^3,J)$ is the infimum of the total (Tanaka--Webster) scalar curvature over all contact forms of volume one.  For CR manifolds $(M^3,J)$ with positive CR Yamabe constant and nonnegative Paneitz operator, the first- and third-named authors~\cite{CY} showed that $\int Q^\prime\leq 16\pi^2$ with equality if and only if $(M^3,J)$ is CR equivalent to the standard CR three-sphere.  The main goal of this note is to refine this statement by giving an integral formula for the total $Q^\prime$-curvature in terms of the Green's function of the CR Laplacian:

\begin{theorem}
 \label{thm:main_thm}
 Let $(M^3,J,\theta)$ be a pseudo-Einstein manifold with positive CR Yamabe constant and nonnegative Paneitz operator.  Given any $p\in M$, it holds that
 \begin{equation}
  \label{eqn:integral_equation}
  \int_M Q^\prime = 16\pi^2 - 4\int_{M} G_L^4\lvert A_{11}\rvert_{\hat\theta}^2 - 12\int_M 3\log(G_L) P_4\log(G_L)
 \end{equation}
 where $G_L$ is the Green's function for the CR Laplacian with pole $p$ and $\hat\theta=G_L^2\theta$.  In particular,
 \[ \int_M Q^\prime \leq 16\pi^2 \]
 with equality if and only if $(M^3,J)$ is CR equivalent to the standard CR three sphere.
\end{theorem}

Theorem~\ref{thm:main_thm} is motivated by similar work in conformal geometry: Gursky~\cite{G} used the total $Q$-curvature to characterize the standard four-sphere among all Riemannian manifolds with positive Yamabe constant and Hang--Yang~\cite{HaY} rederived this result by giving an integral formula for the total $Q$-curvature in terms of the Green's function for the conformal Laplacian.

The key technical difficulty in the proof of Theorem~\ref{thm:main_thm} comes from the potential need to consider the $Q^\prime$-curvature of a contact form which is not pseudo-Einstein.  On the one hand, $\log G_L$ need not be CR pluriharmonic, and hence $G_L^2\theta$ need not be pseudo-Einstein; this problem is overcome by adapting ideas from~\cite{CY}.  On the other hand, estimates for $\log G_L$ are usually derived in CR normal coordinates (cf.\ \cite{HY}), but CR normal coordinates need not be specified in terms of a pseudo-Einstein contact form.  We overcome the latter issue by using Moser's contact form, which is necessarily pseudo-Einstein, as a replacement for CR normal coordinates.

Ignoring these technical difficulties, the idea of the proof of Theorem~\ref{thm:main_thm} is to observe that $\hat\theta:=G_L^2\theta$ has vanishing scalar curvature away from the pole, and hence $\hat Q^\prime$ has a particularly simple expression.  Equation~\eqref{eqn:qprime_transformation} relates $Q^\prime$ and $\hat Q^\prime$ in terms of $P^\prime(\log G_L)$ and $P\bigl((\log G_L)^2\bigr)$.  Using normal coordinates, we can compute these latter functions near the pole $p$, at which point~\eqref{eqn:integral_equation} follows from~\eqref{eqn:qprime_transformation} by integration.  As an upshot of this approach, we relate $\log G_L$ and the Green's function for $P^\prime$; we expect this relation to be useful for future studies of the $Q^\prime$-curvature.

This note is organized as follows: In Section~\ref{sec:moser}, we recall necessary facts about Moser's contact form and use it to relate $Q^\prime$ and $\hat Q^\prime$.  In Section~\ref{sec:integral}, we integrate this relation to prove Theorem~\ref{thm:main_thm}.

\section{Moser's contact form and normal coordinates}
\label{sec:moser}

Moser's normal form for a real hypersurface in $\mathbb{C}^{2}$ (see, e.g., \cite{CM}%
) reads%
\begin{equation*}
v=|z|^{2}-E(u,z,\bar{z})
\end{equation*}

\noindent where $(z,w)$ $\in $ $C^{2},$ $w=u+iv,$ and%
\begin{eqnarray*}
E(u,z,\bar{z}) &=&-c_{42}(u)z^{4}\bar{z}^{2}-c_{24}(u)z^{2}\bar{z}^{4} \\
&&-c_{33}(u)z^{3}\bar{z}^{3}+O(7).
\end{eqnarray*}

\noindent Hereafter we use $O(k)$ to denote $O(\rho ^{k})$ where $\rho $ $:=$
$(|z|^{4}+u^{2})^{1/4}.$

Associated to the defining function%
\begin{equation}
r=\frac{1}{2i}(w-\bar{w})-|z|^{2}+E(u,z,\bar{z}),  \label{N0}
\end{equation}

\noindent we have Moser's contact form%
\begin{eqnarray*}
\theta &=&i\partial r \\
&=&\frac{1}{2}dw-i\bar{z}dz+i(E_{z}dz+E_{u}\frac{1}{2}dw)
\end{eqnarray*}

\noindent in which we have used $E_{w}=E_{u}\frac{1}{2}$ and $E_{z}$ $%
:=\partial E/\partial z,$ $E_{u}$ $:=$ $\partial E/\partial u,$ etc.$.$ We
call coordinates $(z,u)$ for real hypersurface $\{r=0\}$ Moser's normal
coordinates. We are going to compute pseudohermitian quantities with respect
to Moser's contact form in Moser's normal coordinates. Compute%
\begin{equation*}
d\theta =ig_{1\bar{1}}dz\wedge d\bar{z}+\theta \wedge \phi
\end{equation*}

\noindent where%
\begin{eqnarray}
g_{1\bar{1}} &=&1-E_{z\bar{z}}-\lambda E_{u\bar{z}}-\bar{\lambda}%
E_{uz}-|\lambda |^{2}E_{\text{uu}},  \label{N0-0} \\
\phi &=&a_{1}dz+a_{\bar{1}}d\bar{z}  \notag
\end{eqnarray}

\noindent in which%
\begin{eqnarray}
\lambda &=&\frac{\bar{z}-E_{z}}{-i+E_{u}}=i\bar{z}-iE_{z}+\bar{z}E_{u}+O(6)
\label{N0-1} \\
a_{1} &=&\frac{-E_{uz}-\lambda E_{\text{uu}}}{i+E_{u}},\text{ }a_{\bar{1}}=%
\overline{(a_{1})}.  \notag
\end{eqnarray}

\noindent The order counting follows the rule that $z,$ $\bar{z}$ are of
order 1 and $u$ is of order 2. Here we have also used the relation between $%
dw$ and $\theta :$%
\begin{equation*}
dw=\frac{2}{1+iE_{u}}(\theta +i(\bar{z}-E_{z})dz).
\end{equation*}

Take a pseudohermitian coframe%
\begin{align}
\theta ^{1} &:=dz-ia^{1}\theta ,  \label{N1} \\
a^{1} &:=g^{1\bar{1}}a_{\bar{1}}  \notag
\end{align}

\noindent where $g^{1\bar{1}}$ $:=$ ($g_{1\bar{1}})^{-1},$ such that%
\begin{equation}
d\theta =ig_{1\bar{1}}\theta ^{1}\wedge \theta ^{\bar{1}}.  \label{N2}
\end{equation}

\noindent The dual frame $Z_{1}$ (such that $\theta (Z_{1})$ $=$ $0,$ $%
\theta ^{1}(Z_{1})$ $=$ $1,$ and $\theta ^{\bar{1}}(Z_{1})$ $=$ $0)$ reads%
\begin{equation}
Z_{1}=\frac{\partial }{\partial z}+\lambda \frac{\partial }{\partial u}=%
\mathring{Z}_{1}+O(5)\frac{\partial }{\partial u}  \label{N2-1}
\end{equation}

\noindent where $\mathring{Z}_{1}$ $:=$ $\partial _{z}+i\bar{z}\partial
_{u}. $

Differentiating $\theta ^{1}$ from (\ref{N1}) gives%
\begin{equation}
d\theta ^{1}=\theta ^{1}\wedge \mathring{\omega}_{1}^{1}+iZ_{\bar{1}%
}(a^{1})\theta \wedge \theta ^{\bar{1}}  \label{N3}
\end{equation}

\noindent by (\ref{N2}), where%
\begin{equation}
\mathring{\omega}_{1}^{1}=a_{\bar{1}}\theta ^{\bar{1}}-iZ_{1}(a^{1})\theta .
\label{N3-1}
\end{equation}

\noindent Differentiating (\ref{N2}) gives no $\theta $ component of%
\begin{eqnarray*}
&&dg_{1\bar{1}}-g_{1\bar{1}}\mathring{\omega}_{1}^{1}-g_{1\bar{1}}\mathring{%
\omega}_{\bar{1}}^{\bar{1}} \\
&=&[Z_{1}g_{1\bar{1}}-g_{1\bar{1}}a_{1}]\theta ^{1}+\text{ conjugate}
\end{eqnarray*}

\noindent where $\mathring{\omega}_{\bar{1}}^{\bar{1}}$ is the complex
conjugate of $\mathring{\omega}_{1}^{1}.$ It follows that the
pseudohermitian connection form $\omega _{1}^{1}$ reads%
\begin{equation}
\omega _{1}^{1}=\mathring{\omega}_{1}^{1}+(g^{1\bar{1}}Z_{1}g_{1\bar{1}%
}-a_{1})\theta ^{1}.  \label{N4}
\end{equation}

\noindent We also reads from (\ref{N3}) that%
\begin{equation}
A_{\bar{1}}^{1}=iZ_{\bar{1}}(a^{1}).  \label{N5}
\end{equation}

Substituting (\ref{N4}) into the structure equation $d\omega _{1}^{1}$ $=$ $%
Rg_{1\bar{1}}\theta ^{1}\wedge \theta ^{\bar{1}}$ mod $\theta ,$ we obtain
the Tanaka-Webster (scalar) curvature%
\begin{eqnarray}
R &=&Z^{\bar{1}}a_{\bar{1}}+Z_{1}a^{1}+Z^{1}a_{1}  \label{N6} \\
&&-Z^{1}Z^{\bar{1}}g_{1\bar{1}}+a^{1}(Z^{\bar{1}}g_{1\bar{1}})-a^{1}a_{1}. 
\notag
\end{eqnarray}

\noindent where we have used $g^{1\bar{1}}$ to raise the indices, e.g., $Z^{%
\bar{1}}$ $:=$ $g^{1\bar{1}}Z_{1}$ $=$ $\overline{(g^{1\bar{1}}Z_{\bar{1}})}$
$=$ $\overline{(Z^{1})},$ $a^{1}$ $:=$ $g^{1\bar{1}}a_{\bar{1}}.$ We then
compute the lowest order terms of $Z^{\bar{1}}a_{\bar{1}},$ $Z^{\bar{1}}g_{1%
\bar{1}}$ as follows:%
\begin{eqnarray}
Z^{\bar{1}}a_{\bar{1}} &=&-E_{\text{uu}}-iE_{u\bar{z}z}-zE_{\text{uu}z}+\bar{%
z}E_{\text{uu}\bar{z}}-i|z|^{2}E_{\text{uuu}}+O(3),  \label{N7} \\
Z^{\bar{1}}g_{1\bar{1}} &=&-E_{zz\bar{z}}-2i\bar{z}E_{z\bar{z}u}+izE_{uzz}+%
\bar{z}^{2}E_{\text{uu}\bar{z}}+iE_{uz}  \notag \\
&&-2|z|^{2}E_{\text{uu}z}-\bar{z}E_{\text{uu}}-i\bar{z}|z|^{2}E_{\text{uuu}%
}+O(4).  \notag
\end{eqnarray}

\noindent Here we have counted $z,$ $\bar{z}$ of order 1, $u$ of order 2,
and used $g^{1\bar{1}}$ $=$ $1+O(4),$ $\lambda $ $=$ $i\bar{z}-iE_{z}+\bar{z}%
E_{u}$ $+$ $h.o.t.,$ $a_{1}$ $=$ $iE_{uz}-\bar{z}E_{\text{uu}}$ $+$ $h.o.t.$%
, $Z_{1}$ = $\partial _{z}+i\bar{z}\partial _{u}+h.o.t..$ From (\ref{N5}) we
compute%
\begin{equation}
A_{\bar{1}}^{1}=E_{u\bar{z}\bar{z}}-2izE_{\text{uu}\bar{z}}+z^{2}E_{\text{uuu%
}}+O(3).  \label{N7-1}
\end{equation}

\noindent By (\ref{N7}) and alike formulas, we can compute $R$ through (\ref%
{N6}):%
\begin{eqnarray}
R &=&-2E_{\text{uu}}+E_{zz\bar{z}\bar{z}}-2izE_{zz\bar{z}u}+2i\bar{z}E_{z%
\bar{z}\bar{z}u}  \label{N7-2} \\
&&+4|z|^{2}E_{z\bar{z}\text{uu}}-z^{2}E_{\text{uu}zz}-\bar{z}^{2}E_{\text{uu}%
\bar{z}\bar{z}}+2i\bar{z}|z|^{2}E_{\text{uuu}\bar{z}}  \notag \\
&&-2iz|z|^{2}E_{\text{uuu}z}+|z|^{4}E_{\text{uuuu}}+O(3)  \notag
\end{eqnarray}

We can then compute $R_{,1}$ $=$ $Z_{1}R$, $A_{1,\bar{1}}^{\bar{1}}$, and
obtain the pseudo-Einstein tensor as follows:%
\begin{eqnarray}
R_{,1}-iA_{1,\bar{1}}^{\bar{1}} &=&E_{zzz\bar{z}\bar{z}}-4i\bar{z}E_{\text{%
uuu}}+3i\bar{z}E_{zz\bar{z}\bar{z}u}  \label{N8} \\
&&-3iE_{zz\bar{z}u}-2izE_{zzz\bar{z}u}+6|z|^{2}E_{zz\bar{z}\text{uu}}  \notag
\\
&&-3\bar{z}^{2}E_{\text{uu}z\bar{z}\bar{z}}+6\bar{z}E_{z\bar{z}\text{uu}}+6i%
\bar{z}|z|^{2}E_{z\bar{z}\text{uuu}}  \notag \\
&&-3zE_{\text{uu}zz}-z^{2}E_{\text{uu}zzz}-3iz|z|^{2}E_{\text{uuu}zz}  \notag
\\
&&-i\bar{z}^{3}E_{\text{uuu}\bar{z}\bar{z}}+i\bar{z}^{2}E_{\text{uuu}\bar{z}%
}-2\bar{z}^{2}|z|^{2}E_{\text{uuuu}\bar{z}}  \notag \\
&&-6i|z|^{2}E_{\text{uuu}z}+3|z|^{4}E_{\text{uuuu}z}+\bar{z}|z|^{2}E_{\text{%
uuuu}}  \notag \\
&&+i\bar{z}|z|^{4}E_{\text{uuuuu}}+O(2).  \notag
\end{eqnarray}

\noindent From (\ref{N8}) along the $u$-curve (a chain) where $z=0,$ we
conclude that $R_{,1}-iA_{1,\bar{1}}^{\bar{1}}$ $=$ $0$ (terms in $O(2)$ all
vanish because of special structure of Moser's normal form) and does not
vanish identically in general. The reason is that the coefficient of $z$ in $%
E_{zzz\bar{z}\bar{z}}$ is $c_{42}(u)$ which is proportional to the Cartan
tensor.

\bigskip

In general a pseudo-Einstein contact form may not be a ``normalized'' contact
form that gives CR normal coordinates. So we take the contact form
associated to the solution $\psi $ to the complex Monge--Amp\`ere equation:%
\begin{equation}
J[\psi ]:=\det \left[ 
\begin{array}{ccc}
\psi & \psi _{\bar{z}} & \psi _{\bar{w}} \\ 
\psi _{z} & \psi _{z\bar{z}} & \psi _{z\bar{w}} \\ 
\psi _{w} & \psi _{w\bar{z}} & \psi _{w\bar{w}}%
\end{array}%
\right] =1  \label{N8-1}
\end{equation}

\noindent in $\bar{\Omega}$ and $\psi $ $=$ $0$ on $\partial \Omega .$ The
contact form $\theta $ :$=$ $i\partial \psi $ is pseudo-Einstein. We want to
compute $\Delta _{b},$ $P^{\prime },$ $P$ w.r.t. this $\theta ,$ but in
Moser's normal coordinates $(z,u).$ For $r$ having a form of (\ref{N0})
multiplied by $4^{1/3}$, we have%
\begin{equation}
J[r]=1+O(\rho ^{4}).  \label{N9}
\end{equation}

\noindent Lee-Melrose's asymptotic expansion reads%
\begin{equation}
\psi \sim r\sum_{k\geq 0}\eta _{k}(r^{3}\log r)^{k}\text{ near }\partial
\Omega =\{r=0\}\subset C^{2}  \label{N10}
\end{equation}

\noindent with $\eta _{k}$ $\in $ $C^{\infty }(\bar{\Omega}).$ This means
that for large $N,$ $\psi -r\sum_{k=0}^{N}\eta _{k}(r^{3}\log r)^{k}$ has
many continuous derivatives on $\bar{\Omega}$ and vanishes to high order at $%
\partial \Omega .$ It follows from (\ref{N8-1}), (\ref{N9}), and (\ref{N10})
that 
\begin{eqnarray*}
J[r\eta _{0}] &=&1+O(\rho ^{4})\text{ and} \\
\eta _{0} &=&1+O(\rho ^{4}).
\end{eqnarray*}

\noindent So we have%
\begin{align*}
\psi & \sim r\eta _{0}+\eta _{1}r^{4}\log r+h.o.t. \\
& \sim r+O(\rho ^{6}).
\end{align*}

\noindent Similar argument as for $r$ before works for $\psi$. Therefore,
with respect to the pseudo-Einstein contact form defined by $\psi$, we still have%
\begin{eqnarray}
\theta &=&(1+O(\rho ^{4}))\mathring{\theta}+O(\rho ^{5})dz+O(\rho ^{5})d\bar{%
z},  \label{N10-1} \\
\theta ^{1} &=&O(\rho ^{3})\mathring{\theta}+(1+O(\rho ^{8}))dz+O(\rho ^{8})d%
\bar{z},  \notag \\
Z_{1} &=&\mathring{Z}_{1}+O(\rho ^{5})\frac{\partial }{\partial u},  \notag
\\
\omega _{1}^{1} &=&O(\rho ^{2})\mathring{\theta}+O(\rho ^{3})dz+O(\rho ^{7})d%
\bar{z},  \notag \\
A_{\bar{1}}^{1} &=&O(\rho ^{2}),R=O(\rho ^{2}),  \notag \\
g_{1\bar{1}} &=&1+O(\rho ^{4}),\text{ }g^{1\bar{1}}=1+O(\rho ^{4})  \notag
\end{eqnarray}

\noindent in view of (\ref{N0-1}), (\ref{N2-1}), (\ref{N3-1}), (\ref{N4}), (%
\ref{N7-1}), (\ref{N7-2}), and (\ref{N0-0}). Now let $L$ denote the $CR$
Laplacian:%
\begin{equation*}
L:=-4\bigtriangleup _{b}+R
\end{equation*}

\noindent where $\bigtriangleup _{b}$ is the (positive) sublaplacian given by%
\begin{equation}
\bigtriangleup _{b}=Z^{1}Z_{1}-\omega _{1}^{1}(Z^{1})Z_{1}+\text{conjugate.}
\label{N10-2}
\end{equation}

\noindent Let $G_{L}$ denote the Green's function of $L,$ i.e.,%
\begin{equation}
 \label{N11}
 LG_L = -4\bigtriangleup_b G_L + RG_L = 16\delta_p .
\end{equation}
%

\noindent Let $\mathring{P}^{\prime }$ $:=$ $4\mathring{\bigtriangleup}%
_{b}^{2},$ $\mathring{L}$ $:=$ $-4\mathring{\bigtriangleup}_{b}$ denote the $%
P^{\prime }$ operator, the CR Laplacian for the Heisenberg group $\mathbb{H}^1$, 
respectively.  Observe that (cf.\ \cite{BFM})%
\begin{equation*}
\mathring{P}^{\prime }(\log G_{\mathring{L}})=\mathring{P}^{\prime }(\log 
\frac{1}{2\pi \rho ^{2}})=8\pi^2S_{p}
\end{equation*}
where $S_p=S(p,\cdot)$ for $S(p,\cdot)$ the kernel of the orthogonal projection 
$\pi\colon L^2(\mathbb{H}^1)\to\mathcal{P}(\mathbb{H}^1)$ onto the space of CR 
pluriharmonic functions.

\noindent where we have used $G_{\mathring{L}}$ $=$ $\frac{1}{2\pi \rho ^{2}}%
.$ From (\ref{N10-1}) and (\ref{N10-2}) we obtain%
\begin{eqnarray}
\bigtriangleup _{b} &=&(1+O(\rho ^{4}))\mathring{\bigtriangleup}_{b}+O(\rho
^{10})\frac{\partial ^{2}}{\partial u^{2}}+O(\rho ^{4})\frac{\partial }{%
\partial u}  \label{N12} \\
&&+O(\rho ^{5})\frac{\partial }{\partial u}\circ \mathring{Z}_{1}+O(\rho
^{5})\mathring{Z}_{1}\circ \frac{\partial }{\partial u}+O(\rho ^{7})%
\mathring{Z}_{1}  \notag \\
&&+O(\rho ^{5})\frac{\partial }{\partial u}\circ \mathring{Z}_{\bar{1}%
}+O(\rho ^{5})\mathring{Z}_{\bar{1}}\circ \frac{\partial }{\partial u}%
+O(\rho ^{7})\mathring{Z}_{\bar{1}}.  \notag
\end{eqnarray}%
Write%
\begin{equation*}
G_{L}=\frac{1}{2\pi \rho ^{2}}+\omega .
\end{equation*}

\noindent From (\ref{N11}), (\ref{N12}), and (\ref{N10-1}) we obtain%
\begin{equation*}
L\omega =\text{a bounded function near }p.
\end{equation*}

\noindent Therefore from subelliptic regularity theory of $L$, we see that $\omega$ is in
the Folland--Stein space $S^{2,q}$ for any $q>1$, and hence $w\in C^{1,\gamma}$.
In fact, $\omega $ is $C^{\infty }$ smooth~\cite{HY}. Recall that%
\begin{eqnarray}
P^{\prime } &=&4\bigtriangleup _{b}^{2}-8Im\nabla ^{1}(A_{1}^{\bar{1}%
}\nabla _{\bar{1}})-4Re\nabla ^{1}(R\nabla _{1})  \label{N13} \\
&=&\mathring{P}^{\prime }+4(\bigtriangleup _{b}^{2}-\mathring{\bigtriangleup}%
_{b}^{2})  \notag \\
&&-8Im\nabla ^{1}(A_{1}^{\bar{1}}\nabla _{\bar{1}})-4Re\nabla
^{1}(R\nabla _{1}).  \notag
\end{eqnarray}%
\noindent Write 
\begin{eqnarray}
\log G_{L} &=&\log (\frac{1}{2\pi \rho ^{2}}+\omega )  \label{N13-1} \\
&=&\log (\frac{1}{2\pi \rho ^{2}})+\log (1+2\pi \rho ^{2}\omega ).  \notag
\end{eqnarray}%
\noindent We can now compute 
\begin{eqnarray*}
P^{\prime }(\log G_{L}) &=&\mathring{P}^{\prime }(\log (\frac{1}{2\pi \rho
^{2}}))+(P^{\prime }-\mathring{P}^{\prime })(\log (\frac{1}{2\pi \rho ^{2}}))
\\
&&+P^{\prime }(\log (1+2\pi \rho ^{2}\omega )) \\
&=&8\pi^2 S_{p}+\{4(\bigtriangleup _{b}^{2}-\mathring{\bigtriangleup}%
_{b}^{2})-8Im\nabla ^{1}(A_{1}^{\bar{1}}\nabla _{\bar{1}}) \\
&&-4Re\nabla ^{1}(R\nabla _{1})\}(\log (\frac{1}{2\pi \rho ^{2}})) \\
&&+P^{\prime }(\log (1+2\pi \rho ^{2}\omega )).
\end{eqnarray*}

\noindent Since $\omega $ is $C^{\infty }$ smooth, the third term is a
bounded function near $p.$ The second term is also bounded near $p$ in view
of (\ref{N10-1}) and (\ref{N12}). So we conclude that 
\begin{equation}
P^{\prime }(\log G_{L})=8\pi^2S_{p}+\text{a bounded function.}
\label{N13-1-1}
\end{equation}

\noindent Similarly we can show%
\begin{equation}
P((\log G_{L})^{2})=8\pi^2\left(\delta _{p}-S_p\right)+\text{a bounded function.}
\label{N13-1-2}
\end{equation}%
\noindent On the other hand, we reduce computing the most singular term in $%
P_{3}(\log G_{L})$ to computing $P_{3}(\log (\frac{1}{2\pi \rho ^{2}}))$ by (%
\ref{N13-1}). In view of (\ref{N10-1}) we find that the most singular term
in $P_{3}(\log (\frac{1}{2\pi \rho ^{2}}))$ is a constant multiple of $%
\mathring{P}_{3}(\log \rho )$ where $\mathring{P}_{3}$ $=$ $\mathring{Z}_{1}%
\mathring{Z}_{1}\mathring{Z}_{\bar{1}}$ is the $P_{3}$-operator w.r.t. the
Heisenberg group $H_{1}$. Observe that $|z|^{2}-iu$ is a CR function on $%
H_{1}$, i.e.,%
\begin{equation*}
\mathring{Z}_{\bar{1}}(|z|^{2}-iu)=(\partial _{\bar{z}}-iz\partial
_{u})(|z|^{2}-iu)=0.
\end{equation*}%
It follows that the real part of $\log (|z|^{2}-iu)$ is %
CR pluriharmonic. By \cite{Lee} we have%
\[ \mathring{P}_3\left((\log \left| \lvert z\rvert^2 - iu \right|\right) = \mathring{Z}_1\mathring{Z}_1\mathring{Z}_{\bar{1}}\left(\log\left|\lvert z\rvert^2-iu\right|\right) = 0 . \]

\noindent Since $\log\left|\lvert z\rvert^2 - iu\right| = 2\log \rho$, we conclude that
\begin{equation*}
\mathring{P}_{3}(\log \rho )=0.
\end{equation*}

\noindent It follows that%
\begin{eqnarray}
P_{3}(\log G_{L}) &=&\mathring{P}_{3}(\log (\frac{1}{2\pi \rho ^{2}}))
\label{N13-2} \\
&&+(P_{3}-\mathring{P}_{3})(\log (\frac{1}{2\pi \rho ^{2}}))  \notag \\
&&+P_{3}(\log (1+2\pi \rho ^{2}\omega ))  \notag \\
&=&0+(P_{3}-\mathring{P}_{3})(\log (\frac{1}{2\pi \rho ^{2}}))  \notag \\
&&+P_{3}(\log (1+2\pi \rho ^{2}\omega ))  \notag \\
&=&O(\rho ).  \notag
\end{eqnarray}

\noindent by (\ref{N10-1}). So ($\log G_{L})P(\log G_{L})$ has blow-up rate
as $\log \rho $ near the pole $p.$ Hence it is integrable with respect to
the volume $\theta \wedge d\theta $ which has vanishing order $\rho ^{3}$
near $p.$

\section{A formula for the integral of $Q^{\prime }$ curvature}
\label{sec:integral}

Let $\theta$ be a pseudo-Einstein contact form on $(M^{3}, J)$. By~\cite[Proposition~6.1]{CY},
for any $\Upsilon\in C^\infty(M)$, it holds that $\hat{\theta}:=e^\Upsilon\theta$ satisfies
\begin{eqnarray}
e^{2\Upsilon}\hat{Q}^{\prime} &=&Q^{\Upsilon}+P^{\prime }(\Upsilon)+\frac{1}{2}%
P(\Upsilon^{2})  \label{Q1} \\
&&-\Upsilon P(\Upsilon)-16Re(\nabla ^{1}\Upsilon)(P_{3}\Upsilon)_{1} 
\notag
\end{eqnarray}

\noindent where $P_{3}$ is the operator characterizing $CR$ pluriharmonics.
Recall that $P(\Upsilon)$ $=$ $4\nabla ^{1}(P_{3}\Upsilon)_{1}.$

Let $G_{L}$ be the Green's function of the $CR$ Laplacian (we assume $Y(J)$ $%
>$ $0).$ Set $\hat{\theta}$ $=$ $G_{L}^{2}\theta .$ Then $\hat{\theta}$ has
vanishing scalar curvature away from the pole $p.$ In particular, we have%
\begin{equation*}
\hat{Q}^{\prime }=-4|\hat{A}_{11}|_{\hat{\theta}}^{2}
\end{equation*}

\noindent away from the pole $p.$ Plugging this into (\ref{Q1}), we see that
away from $p,$%
\begin{eqnarray}
-4G_{L}^{4}|\hat{A}_{11}|_{\hat{\theta}}^{2} &=&Q^{\prime }+2P^{\prime
}(\log G_{L})  \label{Q2} \\
&&+2P((\log G_{L})^{2})-4(\log G_{L})P(\log G_{L})  \notag \\
&&-64Re(\nabla ^{1}\log G_{L})(P_{3}(\log G_{L}))_{1}.  \notag
\end{eqnarray}

Now assume $(M^{3},$ $J)$ is embedded in $\mathbb{C}^{2}.$ Take $\theta $ to be the
pseudo-Einstein contact form associated to the solution to complex
Monge--Amp\`ere equation (\ref{N8-1}). We look at the order of $G_{L}^{4}|\hat{A%
}_{11}|_{\hat{\theta}}^{2}$ near $p.$ The transformation law of torsion reads%
\begin{equation}
\hat{A}_{11}=G_{L}^{-2}(A_{11}+2i(\log G_{L})_{,11}-4i(\log G_{L})_{,1}(\log
G_{L})_{,1}  \label{Q2-1}
\end{equation}%
\noindent (see~\cite[p.\ 421]{Lee}). Recall $\mathring{Z}_{1}$ $:=$ $\partial
_{z}+i\bar{z}\partial _{u}.$ Observe that 
\begin{eqnarray*}
\mathring{Z}_{1}\log \rho ^{4} &=&\frac{2\bar{z}}{|z|^{2}-iu}, \\
\mathring{Z}_{1}\mathring{Z}_{1}\log \rho ^{4} &=&\frac{-4\bar{z}^{2}}{%
(|z|^{2}-iu)^{2}}=-(\mathring{Z}_{1}\log \rho ^{4})^{2}.
\end{eqnarray*}%
\noindent Therefore we have%
\begin{equation}
\mathring{Z}_{1}\mathring{Z}_{1}\log \frac{1}{2\pi \rho ^{2}}-2(\mathring{Z}%
_{1}\log \frac{1}{2\pi \rho ^{2}})^{2}=0  \label{Q3}
\end{equation}%
\noindent It follows from (\ref{N10-1}) and (\ref{Q3}) that%
\begin{eqnarray}
A_{11} &=&O(\rho ^{2})  \label{Q4} \\
2i(\log G_{L})_{,11}-4i(\log G_{L})_{,1}(\log G_{L})_{,1} &=&O(\rho ^{2}) 
\notag
\end{eqnarray}%
\noindent near $p.$ So from (\ref{Q2-1}) and (\ref{Q4}), we learn that%
\begin{equation}
G_{L}^{4}|\hat{A}_{11}|_{\hat{\theta}}^{2}=O(\rho ^{4})  \label{Q5}
\end{equation}%
\noindent near $p.$ By (\ref{N13-2}), we obtain that the last two terms in (%
\ref{Q2}) are $L^{1}$ and bounded near $p$, respectively.  In view of (\ref%
{N13-1-1}), (\ref{N13-1-2}), (\ref{Q5}), and (\ref{Q2}), we then have%
\begin{eqnarray}
&&2P^{\prime }(\log G_{L})+2P((\log G_{L})^{2})  \label{Q6} \\
&=&16\pi^2\delta _{p}-Q^{\prime }-4G_{L}^{4}|\hat{A}_{11}|_{\hat{%
\theta}}^{2}  \notag \\
&&+4(\log G_{L})P(\log G_{L})+64Re(\nabla ^{1}\log G_{L})(P_{3}(\log
G_{L}))_{1}.  \notag
\end{eqnarray}

\noindent in the distribution sense. Integrating the last term in (\ref{Q6})
gives%
\begin{equation}
-16\int (\log G_{L})P(\log G_{L})+64Re\oint_{\text{around }%
p}(\log G_{L})P_{3}(\log G_{L})i\theta \wedge \theta ^{1}.  \label{N14}
\end{equation}

\noindent Here we have omitted the lower index ``1'' for the $P_{3}$ term. The
boundary term in (\ref{N14}) vanishes by (\ref{N13-2}) and that $\theta
\wedge \theta ^{1}$ has vanishing order of $\rho ^{3}$ near $p.$ Applying (%
\ref{Q6}) to the constant function $1$ yields%
\begin{equation*}
0=16\pi ^{2}-\int Q^{\prime }-4\int G_{L}^{4}|\hat{A}_{11}|_{\hat{\theta}%
}^{2}-12\int (\log G_{L})P(\log G_{L})
\end{equation*}

\noindent by (\ref{N14}). Here notice that in the distribution sense, 
\begin{equation*}
2P^{\prime }(\log G_{L})(1)=2\int (\log G_{L})P^{\prime }(1)=0
\end{equation*}%
\noindent since $P^{\prime }(1)$ $=$ $0.$ Similarly we get $2P((\log
G_{L})^{2})(1)$ $=$ $0$ since $P(1)$ $=$ $0.$ Assuming $P\geq 0,$ we get that%
\begin{eqnarray*}
\int Q^{\prime } &=&16\pi ^{2}-4\int G_{L}^{4}|\hat{A}_{11}|_{\hat{\theta}%
}^{2}-12\int (\log G_{L})P(\log G_{L}) \\
&\leq &16\pi ^{2}.
\end{eqnarray*}
Moreover, equality holds if and only if $\hat{A}_{11}\equiv0$ and $\log G_L$ is pluriharmonic. 
Since also $\hat R\equiv0$, we conclude that $(M\setminus\{p\},\hat\theta)$ is isometric to 
the Heisenberg group $\mathbb{H}^1$.  Indeed, the developing map identifies the universal cover 
of $M\setminus\{p\}$ with $\mathbb{H}^1$, while the fact that a neighorhood of $p$ (equivalently, 
a neighborhood of infinity in $(M\setminus\{p\},\hat\theta)$) is simply connected implies that 
the covering map is trivial.  By adding back the point $p$, we conclude that $(M,J)$ is CR 
equivalent to the standard CR three-sphere.

\end{document}